\documentclass[12pt]{article}

\usepackage{amsmath}
\usepackage{amssymb}
\usepackage{amscd}
\usepackage[all]{xy}

\newtheorem{theorem}{Theorem}[section]

\newtheorem{proposition}[theorem]{Proposition}
\newtheorem{corollary}[theorem]{Corollary}
\newtheorem{lemma}[theorem]{Lemma}
\newtheorem{remark}[theorem]{Remark}

\newenvironment{proof}{\begin{trivlist} \item[]{\bf Proof.}}
{\par\hfill $\square$\end{trivlist}}

\newcommand{\dist}{{\rm dist}}
\renewcommand{\L}{{\cal L}}
\newcommand{\Lone}{{{\rm L}^1}}
\newcommand{\Loneloc}{{\rm L}^1_{\rm loc}}
\newcommand{\Ltwo}{{{\rm L}^2}}

\newcommand{\Lp}{{\rm L}^p}
\newcommand{\Linfty}{{{\rm L}^\infty}}
\newcommand{\Lip}{{\rm Lip}}
\newcommand{\C}{\mathbb{C}}

\renewcommand{\P}{\mathbb{P}}
\newcommand{\R}{\mathbb{R}}

\newcommand{\ddc}{{{\rm dd}^{\rm c}}}
\renewcommand{\d}{{\rm d}}

\newcommand{\Wcap}{{\mbox{\rm W-cap}}}

\newcommand{\DSH}{{\rm DSH}}

\title{Decay of correlations and central limit theorem for meromorphic maps}
\author{Tien-Cuong Dinh and Nessim Sibony}
\date{October 1, 2004}

\begin{document}

\maketitle

\begin{abstract} Let $f$ be a dominating meromorphic self-map 
of large topological degree 
on a compact K\"ahler 
manifold. We give a new construction 
of the equilibrium measure $\mu$ of $f$ and prove that $\mu$ is 
exponentially mixing. Then, we deduce the central limit theorem for 
Lipschitzian observables.
\end{abstract}

\section{Introduction}

Let $(X,\omega)$ be a compact K\"ahler manifold of dimension $k$. 
The K\"ahler form $\omega$ is normalized by $\int\omega^k=1$.
Consider a 
dominating meromorphic self-map $f$ on $X$, i.e. a map whose image contains 
a Zariski open set. 
We assume that the topological degree
$d_t$ of $f$ is strictly larger than the dynamical degree $d_{k-1}$ of order $k-1$.
For such a map, one knows that $d_t^{-n}f^{n*}(\omega^k)$ converges to an 
invariant measure $\mu$. 
The measure $\mu$ is the unique measure of maximal entropy and does not charge
pluripolar sets.
Two different constructions of $\mu$ are given in 
\cite{Guedj,DinhSibony2}. See also \cite{Sibony,RussakovskiiShiffman, BriendDuval} 
for the case of projective spaces, \cite{DinhSibony1} for polynomial-like maps
and \cite{Gromov, Yomdin, DinhSibony3, DinhSibony4} for the bounds of the entropy.

The $\ddc$-method that we used in different situations 
\cite{DinhSibony1, DinhSibony2, DinhSibony5, DinhSibony6} 
leads us to precise controls of the 
convergence. For the maps that we consider, it allowed to prove that $\mu$ 
is exponentially mixing for q.p.s.h. observables; in particular, for ${\cal C}^2$
observables. We then deduce the central limit theorem (CLT for short) 
via a classical result of
Gordin-Liverani. See \cite{FornaessSibony} for 
holomorphic endomorphisms of $\P^k$, \cite{Denker, Haydn} and the references 
therein for the case of dimension 1.

In this article, we develop another method (d-method) which gives a new construction
of the measure $\mu$ and allows to prove analogous statistical results for 
Lipschitzian observables. We obtain in particular the following result.

\
\\
{\it The measure $\mu$ is exponentially 
mixing in the following sense: for every $\epsilon>0$, there exists a
constant $A_\epsilon>0$ such that
$$|\langle\mu,(\psi\circ f^n) \varphi\rangle - \langle \mu,\psi\rangle
\langle \mu,\varphi\rangle |\leq A_\epsilon 
\|\psi\|_\infty\|\varphi\|_\Lip (d_{k-1}+\epsilon)^{n/2} 
d_t^{-n/2}$$
for all $n\geq 0$, every function $\psi\in\Linfty(\mu)$ and 
every Lipschitzian function 
$\varphi$. In particular, if a real-valued Lipschitzian function $\varphi$ is not a
coboundary and verifies $\langle\mu,\varphi\rangle=0$ then it satisfies the 
CLT theorem.}

\
\\
This result is in fact valid for a larger space of test functions
(Theorem 6.1 and Corollary 6.2).
We recall in 
Section 2 the Gordin-Liverani theorem, in Section 3 some properties of the
Sobolev space $W^{1,2}$. The space of test functions $W^{1,2}_*$ and a subspace
$W_{**}^{1,2}$ are introduced in Section 4. 
This is the key point of the method. 
The new space $W_*^{1,2}$ seems to be useful for complex dynamics in
 several variables. In complex dimension 1 it coincides with the space $W^{1,2}$.
 In higher dimension it takes into account, the fact that
  not all currents of bidegree $(1,1)$ are closed. It enjoys the composition
 properties under meromorphic maps, useful for a space of observables.
In Section 5, we give the new 
construction of $\mu$ and in Section 6 its statistical properties. Note that 
the geometric decay of correlations for H\'enon maps was recently proved in 
\cite{Dinh} using the $\ddc$-method.

\
\\
{\bf Notations.} We will use different subspaces of $\Lone(X)$. Most of them carry
a canonical (quasi)-norm. 
For the space ${\cal C}^0(X)$ of continuous functions we use 
the sup-norm, and for the space $\Lip(X)$ of Lipschitzian functions 
the norm $\|.\|_\Lone +\|.\|_\Lip$ where
$\|\varphi\|_\Lip:=\sup_{x\not=y}|f(x)-f(y)|\dist(x,y)^{-1}$.
We use the norm $\|.\|_E+\|.\|_F$ for
the intersection $E\cap F$ of $E$ and $F$
and write ${\rm L}^p$ instead of ${\rm L}^p(X)$ when there is no confusion.
The $\Lp$ norm of a form is the sum of $\Lp$ norms of its coefficients for a
fixed atlas of $X$. The topology, that we consider, 
is {\bf a weak topology}. That is, $\varphi_n\rightarrow\varphi$ in $E$ if 
$\varphi_n\rightarrow\varphi$ in the sense of distributions 
and $(\|\varphi_n\|_E)$
is bounded. The continuity of linear operators $\Lambda:E\rightarrow F$ is with
respect to these topologies of $E$ and $F$. 
The inclusion map $E\subset \Lone(X)$ 
is always bounded for the associated (quasi)-norms 
and continuous in our sense.

\section{Gordin-Liverani theorem}

Let $(X,\mu)$ be a probability space. Let $f:X\rightarrow X$ be
a measurable map in the sense that $f(A)$ and $f^{-1}(A)$ are
measurable for every measurable set $A$. Assume
that $\mu$ is invariant ($f_*(\mu)=\mu$) and ergodic.
The classical Birkhoff ergodic theorem says that if $\varphi$ is 
in $\Lone(\mu)$ then $\frac{1}{n}\sum_{i=0}^{n-1}
\varphi\circ f^i$ converges almost everywhere to $\langle\mu,\varphi\rangle$.
In particular, when $\varphi$ is a function in
$\Lone(\mu)$ such that
$\langle\mu,\varphi\rangle=0$, then $\frac{1}{n}\sum_{i=0}^{n-1}
\varphi\circ f^i$ converges almost everywhere to $0$.

Assume that $\langle\mu,\varphi\rangle=0$. 
Recall that $\varphi$ is a {\it coboundary} 
if there exists a measurable function 
$\psi$ such that $\varphi=\psi\circ f -\psi$. It is easy to check for a  
coboundary $\varphi$ that 
$\frac{1}{\sqrt{n}} \sum_{i=0}^{n-1}
\varphi\circ f^i$ converges to 0 almost surely.
We say that a real-valued function $\varphi\in\Lone(\mu)$ satisfies
the {\it central limit theorem} if
$\frac{1}{\sqrt{n}} \sum_{i=0}^{n-1}
\varphi\circ f^i$
converges in law to a Gaussian random variable of zero mean and
variance $\sigma>0$. That is, for every interval $I\subset \R$,
\begin{eqnarray*}
\lim_{n\rightarrow\infty} \mu \left\{\frac{1}{\sqrt{n}} \sum_{i=0}^{n-1}
\varphi\circ f^i \in I\right\} =\frac{1}{\sqrt{2\pi}\sigma}\int_I
\exp\left({-\frac{x^2}{2\sigma^2}}\right) \d x.
\end{eqnarray*}
Such a function is not a coboundary.

Observe that $f^*$ acts on $\Ltwo(\mu)$ and preserves the $\Ltwo$ norm.
Hence, we can consider the adjoint $\Lambda$ of $f^*$, sometimes called
the {\it Perron-Frobenius operator}. 
Consider now a real-valued function $\varphi\in\Linfty(\mu)$ such that
$\langle\mu,\varphi\rangle=0$ which is not a coboundary.
The theorem of Gordin-Liverani \cite{Gordin}, 
\cite[pp.59, 67]{Liverani} implies that if 
\begin{eqnarray}
\sum_{n\geq 0} \|\Lambda^n
\varphi\|_{\Lone(\mu)}<+\infty
\end{eqnarray}
then $\varphi$ satisfies the CLT
with
\begin{eqnarray*}
\sigma^2=-\langle\mu,\varphi^2\rangle +2 \sum_{n=0}^\infty \langle\mu,
\varphi(\varphi\circ f^n)\rangle.
\end{eqnarray*}

\section{Sobolev space}

In this section we recall some properties of the Sobolev space $W^{1,2}$ that 
will be used later on, see \cite[pp. 23-27]{Hebey}. 

Consider a compact Riemannian manifold $X$ of dimension $m\geq 2$. 
The Sobolev space $W^{1,2}$ is the space
of real-valued functions $\varphi$ in $\Ltwo(X)$ such that 
$\d\varphi$, which is defined in the sense of currents, has $\Ltwo$ coefficients.

Let $\L$ denote the Lebesgue measure on $X$ that we normalize by $\|\L\|=1$. 
Define for each $\Lone$ function $\varphi$ its mean value 
$m(\varphi):=\int\varphi\d\L$. We have the following Poincar\'e-Sobolev 
inequality: for every real number $p$, $1\leq p\leq 2m/(m-2)$, there exists a
constant $c>0$ such that
$$\left(\int_X|\varphi-m(\varphi)|^p\d\L\right)^{1/p}\leq c\|\d\varphi\|_{\Ltwo}$$
for every $\varphi\in W^{1,2}$. 
In particular, this holds
for $p=1$ or $2$. One deduces that
\begin{eqnarray}
\|\varphi\|_{{\rm L}^p}\leq |m(\varphi)|+c\|\d\varphi\|_{\Ltwo}\leq \|\varphi\|_\Lone
+c\|\d\varphi\|_\Ltwo.
\end{eqnarray}
Hence, $W^{1,2}\subset {\rm L}^p$. The associated inclusion map is bounded.   
\\

We will prove the following proposition for the reader's convenience, see also
\cite{David, HarveyPolking}.

\begin{proposition} Let $I\subset X$ be a compact set whose Hausdorff
$(m-1)$-dimensional measure $H^{m-1}(I)$ is zero. Let $\varphi$ be 
a real-valued function in $\Loneloc(X\setminus I)$. Assume that the coefficients of
$\d\varphi$
are in $\Ltwo(X\setminus I)$. Then $\varphi$ belongs to $W^{1,2}$. 
Moreover, there exist a compact set $M\subset X\setminus I$ and a constant $c>0$, both
independent of $\varphi$, such that
$$\|\varphi\|_{\Lone(X)}\leq c\left(\|\varphi\|_{\Lone(M)} +\|\d\varphi\|_{\Lone(X)}
\right).$$
\end{proposition}
\begin{proof} First, assume that $\varphi\in \Lone(X)$. Then 
$\d\varphi$ defines a flat current
of degree $1$ on $X$ \cite{Federer}. 
Since $\d\varphi$ has coefficients in $\Ltwo(X\setminus I)$,
$\d\varphi_{|X\setminus I}$ defines also a flat current on $X$. It follows that 
$\d\varphi-\d\varphi_{|X\setminus I}$ is a flat current of dimension $m-1$
with support in $I$. Since $H^{m-1}(I)=0$, 
by the support theorem \cite[4.1.20]{Federer}, 
this current vanishes. Hence,
$\d\varphi=\d\varphi_{X\setminus I}$. Using a regularization and the 
Poincar\'e-Sobolev inequality, we can approach $\varphi$
in $\Lone(X)$ by a bounded sequence $(\varphi_n)\subset W^{1,2}$.
It follows that $\varphi\in W^{1,2}$.

Now, it is sufficient to prove the estimate in the proposition. 
Consider a chart $D$ of
$X$ that we identify to $]-2,2[\times B_2$ where $B_r$ is the ball of 
center $0$ and of radius $r$ in $\R^{m-1}$. Consider also the
local coordinates $x=(x_1,x')$ with $x\in]-2,2[$ and $x'\in B_2$. 
We can assume that $I\cap D\subset \{|x_1|<1\}$ since $H^{m-1}(I)=0$. Let 
$\chi$ be a smooth positive function with compact support in $]-2,2[$
such that $\chi=1$ on $[-1,1]$. 
We need only to bound the $\Lone$ norm
of $\chi\varphi$ in $]-2,2[\times B_1$ by a multiple of the $\Lone$ norm
of $\d(\chi\varphi)$ in $]-2,2[\times B_1\setminus I$. 

Consider the projection $\pi:D\rightarrow B_2$ defined by 
$\pi(x):=x'$. Then, $\pi(I\cap D)$ is a closed subset of $B_2$
with  $H^{m-1}(\pi(I\cap D))=0$. Hence, since $\varphi\in \Loneloc(X\setminus I)$, 
for almost every $a\in B_1$, 
the restriction of $\chi\varphi$ and $\d(\chi\varphi)$ to 
$L_a:=\pi^{-1}(a)$ are of class $\Lone$. 
Since $(\chi\varphi)_{|L_a}$ has compact support in $]-2,2[$, 
we have 
$$\|(\chi\varphi)_{|L_{a}}\|_\Lone\leq 
4\|\d(\chi\varphi)_{|L_{a}}\|_\Lone.$$  
The Fubini theorem implies that the $\Lone$ norm of $\chi\varphi$ on $]-2,2[\times
B_1$ is bounded by 4 times the $\Lone$ norm of $\d(\chi\varphi)$
on $]-2,2[\times B_1\setminus I$. 
\end{proof}

\section{A space of test functions}

Let $(X,\omega)$ be a compact K\"ahler manifold of dimension $k$
and $\L:=\omega^k$. We assume that $\|\L\|=1$. The space
of test functions that we consider is a subspace $W_*^{1,2}$ of the Sobolev
space $W^{1,2}$ of $X$. A real-valued function $\varphi\in W^{1,2}$ is in 
$W_*^{1,2}$ if there exists a closed current $\Theta$ of bidegree $(1,1)$ 
on $X$ such that  
$i\partial\varphi\wedge \overline\partial \varphi\leq \Theta$.
The current $\Theta$ is necessarily positive since $i\partial\varphi\wedge \overline
\partial \varphi$ is positive.
We identify two fonctions in $W_*^{1,2}$ if they are equal out of a subset of
$\L$ measure zero.
For $\varphi\in W^{1,2}_*$ define 
$$\|\varphi\|_*:=|m(\varphi)| +\inf\{\|\Theta\|^{1/2},\ \Theta \mbox{ as above}\}.$$   
Observe that $\|.\|_*$ is a quasi-norm. Indeed, the Cauchy-Schwarz inequality
implies that
$$i(\partial\varphi+\partial\psi)\wedge
(\overline\partial\varphi+\overline\partial\psi)
\leq 2i\partial\varphi\wedge \overline\partial\varphi + 2i\partial\psi\wedge
\overline\partial\psi,$$
hence $\|\varphi+\psi\|_*\leq \sqrt{2}(\|\varphi\|_*+\|\psi\|_*)$.

A subset of $W_*^{1,2}$ is bounded if $\|.\|_*$ is bounded on this set.
We have $W_*^{1,2}\subset W^{1,2}\subset \Ltwo$
and the associated inclusion maps are bounded and continuous.
Observe also that real-valued Lipschitzian functions belong to this space and when
$k=1$ we have $W^{1,2}_*=W^{1,2}$.  
\\

We have the following proposition

\begin{proposition} Let $\chi:\R\rightarrow\R$ be a Lipschitzian function and 
$\varphi$ be a function in $W_*^{1,2}$. Then $\chi\circ\varphi$ belongs to
$W_*^{1,2}$. In particular, $\varphi^+:=\max(\varphi,0)$, $\varphi^-:=\max(-\varphi,0)$
and $|\varphi|$ belong to $W_*^{1,2}$. Moreover, there exists a constant $c>0$
independent of $\varphi$ such that
$\|\varphi^+\|_*\leq c\|\varphi\|_*$,  $\|\varphi^-\|_*\leq c\|\varphi\|_*$
and $\||\varphi|\|_*\leq c\|\varphi\|_*$. 
If $\varphi_1$, $\varphi_2$ are in $W_*^{1,2}$,
we have $\max(\varphi_1,\varphi_2)\in W_*^{1,2}$.
\end{proposition}
\begin{proof} Since $\chi$ is Lipschitzian, there exist 
 $a\geq 0$, $b\geq 0$ such that 
$|\chi(x)|\leq a+b|x|$. We can choose $a:=\chi(0)$ and $b:=\|\chi\|_\Lip$.
Then, $|\chi\circ \varphi|\leq a+b|\varphi|$. It follows
that $\chi\circ\varphi$ is in $\Ltwo$. We also have
$$i\partial(\chi\circ\varphi)\wedge \overline \partial(\chi\circ\varphi)
=i(\chi'\circ\varphi)^2 \partial\varphi\wedge\overline\partial\varphi
\leq ib^2 \partial \varphi\wedge\overline\partial\varphi\leq b^2\Theta$$
where $\Theta$ is the positive closed current associated to $\varphi$. Hence
$\chi\circ\varphi\in W_*^{1,2}$.

We can apply this for the functions $\chi(x)=\max(x,0)$, $\max(-x,0)$ or
$|x|$ and for $a=0$, $b=1$. 
We obtain that $\varphi^+$, $\varphi^-$, $|\varphi|$ belong to $W^{1,2}_*$.
The estimates on their norms
follow the above inequalities and the Poincar\'e-Sobolev inequality (2). 
For $\max(\varphi_1,\varphi_2)$, it is sufficient to write 
$\max(\varphi_1,\varphi_2)=\max(\varphi_1-\varphi_2,0)+\varphi_2$.
\end{proof}

Recall that an $\Lone$ function $\varphi:X\rightarrow\R\cup\{-\infty\}$ is {\it 
quasi-plurisubharmonic} (q.p.s.h. for short) if it is upper semi-continuous
and verifies $\ddc\varphi\geq -c\omega$ in the sense of currents for some 
constant $c>0$. Note that a q.p.s.h. function is defined at every point of $X$. 
Such a function belongs to $\Lp$ for every 
real number $p\geq 1$, see \cite{Demailly, Lelong, Sibony}.
A subset $Y$ of $X$ is called {\it pluripolar} if
it is contained in the pole set $\{\varphi=-\infty\}$ 
of a q.p.s.h. function $\varphi$.

Let $\DSH$ denote the space generated by q.p.s.h. functions. 
That is, d.s.h. functions can be written, outside a plupripolar set, 
as $\varphi=\varphi_1-\varphi_2$
with $\varphi_i$ q.p.s.h. 
We identify two d.s.h. functions if they are
equal outside a pluripolar set. If $\varphi$ is d.s.h. then 
$i\partial\overline\partial \varphi=T^+-T^-$ where $T^\pm$ are positive
closed $(1,1)$-currents. Define a norm on $\DSH$ by 
$$\|\varphi\|_\DSH:=|m_\varphi|+\inf\{\|T^+\|,\ T^\pm \mbox{ as above}\}.$$
Observe that since $T^+$ and $T^-$ are cohomologous, they have the same mass.
Recall that the mass of a positive $(p,p)$-current $S$ is given by 
$\|S\|=\langle S,\omega^{k-p}\rangle$.

\begin{lemma} We have $\Lip(X)\subset W_*^{1,2}$ and 
$\DSH\cap \Linfty \subset W^{1,2}_*$. Moreover, the corresponding
inclusion maps are bounded and continuous.
\end{lemma}
\begin{proof} 
The inclusion $\Lip(X)\subset W_*^{1,2}$ and the properties of the 
associated map are clear. Recall that we use the norm $\|.\|_\Lone+\|.\|_\Lip$
for $\Lip(X)$.

Let $\varphi\in \DSH\cap \Linfty$ and $T^\pm$ as above. 
We have
$$2i\partial \varphi\wedge \overline\partial \varphi= 
i\partial\overline\partial \varphi^2
-2\varphi i\partial \overline\partial \varphi 
\leq  i\partial\overline\partial \varphi^2
+2\|\varphi\|_\infty (T^++T^-).$$
The right hand side is a closed current cohomologous to 
$4\|\varphi\|_\infty T^+$. Its mass is equal to $4\|\varphi\|_\infty \|T^+\|$.
Hence $\varphi\in W_*^{1,2}$. We also deduce that
$$\|\varphi\|_*\leq \|\varphi\|_\infty + 
\sqrt{2\|\varphi\|_\infty\|\varphi\|_\DSH}.$$ 
It follows that the inclusion map   $\DSH\cap \Linfty \subset W^{1,2}_*$
is bounded. The continuity follows.
\end{proof}

Let $\nu$ be a positive measure on $X$. We say that $\nu$ is {\it WPC} 
if its restriction 
to smooth real-valued functions can be 
extended to a linear form on $W_*^{1,2}$ which is bounded and 
continuous with respect to the weak topology
we consider. 
Let $\langle \nu,\varphi\rangle_*$ denote the value of this linear form 
on $\varphi\in W_*^{1,2}$. We have $\langle \nu,\varphi\rangle_* = 
\langle \nu,\varphi\rangle$ for $\varphi$ smooth. 

If a continuous function $\varphi$ can be uniformly approximated by a bounded 
sequence in $W^{1,2}_*$ of smooth functions, then 
$\langle \nu,\varphi\rangle_* = 
\langle \nu,\varphi\rangle$.
If $\nu$ is given by a $(k,k)$-form with coefficient in $\Ltwo$ then $\nu$
is WPC because $W_*^{1,2}\subset \Ltwo$ and we can take $\langle\nu,\varphi\rangle_*
=\langle\nu,\varphi\rangle$. 

In general, we don't know if smooth 
functions are dense in $W_*^{1,2}$ nor in $W^{1,2}_*\cap {\cal C}^0(X)$. 
Hence, we don't know if the extension of 
$\mu$ is unique and a priori we don't have $\langle \nu,\varphi\rangle_* = 
\langle \nu,\varphi\rangle$ on $W_*^{1,2}\cap {\cal C}^0(X)$.
The method, used in \cite{DinhSibony4}
to regularize currents, might give the solution to this question. 

The proof of the following proposition is left to the reader. One may
apply the classical method of regularization using the diffeomorphisms of $X$ or
the holomorphic automorphisms when $X$ is homogeneous.

\begin{proposition} Smooth functions
are dense in $\Lip(X)$ for the strong topology. In particular, we have 
$\langle \nu,\varphi\rangle_*=\langle\nu,\varphi\rangle$ for $\nu$ WPC and 
$\varphi$ Lipschitzian.
When $k=1$ or $X$ is homogeneous, smooth functions
are dense in $W_*^{1,2}$ and every function in $W^{1,2}_*\cap {\cal C}^0(X)$
can be uniformly approximated by a bounded sequence in $W_*^{1,2}$ of smooth 
functions.
In these cases, $\langle \nu,\varphi\rangle_*=\langle\nu,\varphi\rangle$ holds 
for $\nu$ WPC and 
$\varphi$ in $W_*^{1,2}\cap {\cal C}^0(X)$.
\end{proposition}

\begin{lemma} Let $\nu$ be a WPC measure and $\varphi$ be a bounded q.p.s.h. 
function. Then
$\langle\nu,\varphi\rangle_*=\langle\nu,\varphi\rangle$.
\end{lemma}
\begin{proof} We can assume that $\varphi$ is strictly negative. By Demailly's 
regularization theorem \cite{Demailly}, 
there exist negative smooth functions $\varphi_n$
decreasing to $\varphi$ such that $i\partial\overline\partial \varphi_n
\geq -c\omega$ where $c>0$ is independent of $n$. It follows that 
$\varphi_n\rightarrow\varphi$ in $\DSH\cap\Linfty(X)$ for the weak topology that
we consider. Lemma 4.2 implies that 
$\varphi_n\rightarrow\varphi$ in $W_*^{1,2}$. Hence
$$\langle\nu,\varphi\rangle = \lim \langle\nu,\varphi_n\rangle
= \lim \langle\nu,\varphi_n\rangle_*=\langle\nu,\varphi\rangle_*.$$
\end{proof}

\begin{proposition} Let $\nu=\alpha +\partial u+\overline\partial v$ 
be a positive measure on $X$
where $\alpha$, $u$ and $v$ are forms
with coefficients in $\Ltwo$. Then $\nu$ is WPC. 
\end{proposition}
\begin{proof} If $\varphi$ is smooth, we have
$$\langle \nu,\varphi\rangle = \langle \alpha,\varphi\rangle + 
\langle u,\partial \varphi\rangle +\langle v,\overline\partial \varphi\rangle.$$
The right hand side can be extended to a continuous linear form on $W^{1,2}_*$.
Hence $\nu$ is WPC. 
\end{proof}

Lemma 4.2 and Proposition 4.5 imply the following corollary.

\begin{corollary} Let $\nu=\alpha+i\partial \overline\partial u\wedge \beta$ 
be a positive measure on $X$
where $\alpha$ and $\beta$ are  closed real-valued continuous forms
of bidegree $(k,k)$ and $(k-1,k-1)$ respectively, 
and $u$ is a bounded d.s.h. function. Then $\nu$ is WPC.
\end{corollary}

Let $Y\subset X$ be a Borel set. Define the {\it W-capacity} of $Y$ by
$$\Wcap(Y):= \sup \{\nu(Y), \ \nu \mbox{ a WPC probability measure}\}.$$
Let $W^{1,2}_{**}$ denote the space of functions $\varphi\in W_*^{1,2}$ which are 
continuous outside a compact set of W-capacity zero.
We associate to this space the quasi-norm $\|.\|_*$ and the corresponding topology
(see Introduction).
We have the following proposition.

\begin{proposition} If $Y$ is a pluripolar subset of $X$, 
then $\Wcap(Y)=0$. 
In particular, countable unions of proper analytic subsets of $X$ have W-capacity zero. 
\end{proposition}
\begin{proof} 
Let $\psi$ be a q.p.s.h. function such that $\psi<-4$, 
$i\partial\overline\partial \psi
\geq -\omega$ and $\psi=-\infty$ on $Y$. 
Let $\chi$ be an increasing convex function on $\R\cup\{-\infty\}$ such that 
$\chi=-3$ on $[-\infty,-4]$, $\chi(x)=x$ on $[-2,+\infty[$ and $0\leq\chi'\leq 1$ on
$[-4,-2]$. Define $\chi_n(x):=\chi(x+n)-n$ and $\psi_n:=\chi_n\circ\psi$. 
The functions $\psi_n$ satisfy $\psi_n\leq -3$
and decrease to $\psi$. We have
$$i\partial\overline\partial\psi_n
=(\chi''_n\circ\psi)i\partial\psi\wedge \overline\partial\psi + (\chi'_n\circ\psi)
i\partial\overline\partial\psi \geq (\chi'_n\circ\psi) i\partial\overline\partial
\psi \geq -\omega.$$

Define $\varphi:=-\log(-\psi)$ and $\varphi_n:=-\log(-\psi_n)$. 
The functions $\varphi_n$ decrease to $\varphi$ and
satisfy $\psi_n<-1$. The function 
$\varphi$ is in $\Ltwo$ since $\psi\leq\varphi<0$.
We have
$$i\partial\varphi\wedge \overline\partial\varphi=i\psi^{-2}\partial\psi
\wedge \overline\partial\psi$$
and
$$i\partial\overline\partial \varphi = -i\psi^{-1}\partial\overline\partial \psi +
i\psi^{-2}\partial\psi\wedge\overline\partial\psi.$$
It follows that $i\partial\varphi\wedge \overline\partial\varphi\leq
i\partial\overline\partial \varphi +i \psi^{-1}\partial\overline\partial 
\psi\leq  
i\partial\overline\partial \varphi +\omega$. The positive closed current 
$i\partial\overline\partial\varphi +\omega$ is cohomologous to $\omega$; its mass 
is equal to $1$. 
The functions $\varphi_n$ satisfy a similar inequality. Hence $\varphi_n$
converge to $\varphi$ in $W^{1,2}_*$. 
This and  Lemma 4.4 imply that
$$\langle \nu,\varphi\rangle =\lim \langle \nu,\varphi_n\rangle
=\lim \langle \nu,\varphi_n\rangle_* = \langle\nu,\varphi\rangle_*$$
for every WPC measure $\nu$. 
Hence $\nu(Y)=0$ and $\Wcap(Y)=0$
since $\varphi=-\infty$ on $Y$. 
When $Y$ is a proper analytic subset of $X$, by \cite{DinhSibony2}, $Y$ is pluripolar.
Then $\Wcap(Y)=0$.
\end{proof}

\begin{remark}\rm
{\bf (a)} 
Proposition 4.7 shows that the integral $\langle\nu,\varphi\rangle$ makes sense for
every WPC measure $\nu$ and every bounded d.s.h. function $\varphi$.
Such a function is defined out of a pluripolar set which has $\nu$ measure zero. 

{\bf (b)} The above proof shows that if $\psi$ is a strictly 
negative q.p.s.h. function,
then the function $\varphi:=-\log(-\psi)$ is in $W_*^{1,2}$.
The function $\varphi$ satisfies $\langle\nu,\varphi\rangle_* = 
\langle\nu,\varphi\rangle$ and
could be unbounded.
\end{remark}

\section{Equilibrium measure}

Let $(X,\omega)$ be a compact K\"ahler manifold of dimension $k$ as in Section 4.
Consider a dominating meromorphic self-map $f$ on $X$. Let 
$I_n$ be the indeterminacy set of $f^n$. Define 
$$d_{p,n}:=\int_{X\setminus I_n} f^{n*}(\omega^p)\wedge \omega^{k-p}.$$
In \cite{DinhSibony3, DinhSibony4}, we proved that the limit 
$d_p:=\lim (d_{p,n})^{1/n}$ always exists. It is called the {\it dynamical degree}
of order $p$ of $f$ and is a bimeromorphic invariant. The last 
degree $d_t:=d_k$ is the {\it topological degree} of $f$. Define 
$\delta_n:=d_{k-1,n}$. For example, if $f$ is a meromorphic map in $\P^2$,
$\delta_n$ is just the algebraic degree of $f^n$. 
If $f$ is holomorphic in $\P^k$ of algebraic degree $d$, then $\delta_n =d^{(k-1)n}$,
$d_{k-1}=d^{k-1}$ and $d_t=d^k$.

Let $M$ be the family of positive measures on $X$ 
which does not charge proper analytic sets. If $\nu$ is such a measure,
$f^{n*}(\nu)$ is well defined and $\|f^{n*}(\nu)\|=d_t^n\|\nu\|$
since $f^n$ is locally biholomorphic outside an analytic set and a generic fiber
of $f^n$ contains exactly $d_t^n$ points.

Let $M(\alpha)$, $\alpha>0$, denote
the set of probability measures $\nu\in M$
such that
$|\langle \nu, \varphi \rangle| \leq \alpha\|\varphi\|_*$ for every 
$\varphi\in W_*^{1,2}\cap \Linfty(X)$ which is 
continuous out of an analytic set. 
By Poincar\'e-Sobolev inequality (2), 
we have $\|\varphi\|_\Ltwo\leq c\|\varphi\|_*$, $c>0$. Hence,
if $\|\nu\|_\Ltwo\leq \alpha c^{-1}$ then $\nu\in M(\alpha)$.
\\

The main result of this section is the following theorem.

\begin{theorem} Let $X$ and $f$ be as above. Assume that $d_t>d_{k-1}$. 
Let $\nu_n$ be a probability measure in $M(\alpha_n)$. If 
$\lim \alpha_n^2\delta_nd_t^{-n}=0$
then $d_t^{-n} f^{n*}(\nu_n)$ converge weakly to a measure $\mu$. Moreover, $\mu$
does not depend on $(\nu_n)$ and $\mu$ is WPC.  
\end{theorem}
\noindent
Of course, we have $\mu=\lim d_t^{-n} f^{n*}(\omega^k)$. This is the equilibrium 
measure of $f$.
Theorem 5.1 and Proposition 4.7 show that $\mu$ does not charge pluripolar sets.
The previous relation implies that $f^*(\mu)=d_t\mu$ and $f_*(\mu)=\mu$.
\\

Consider a Zariski open set  $\Omega_n$ of $X$ 
such that $f$ is locally  biholomorphic and proper on $f^{-1}(\Omega_n)$. 
Let $\Theta$ be a positive closed current of bidegree $(p,p)$ on $X$. 
Then we can define $\Theta_n:=(f^n)_*(\Theta)$ on $\Omega_n$. 
If this current has finite mass, its trivial extension
is a positive closed current that we denote by $(f^n)_\star(\Theta)$.
The following lemma shows that this is the case, see also 
\cite{DinhSibony3, DinhSibony4}.
The choice of $\Omega_n$ is not important.

\begin{lemma} There exists a constant $c>0$ which depends only on $(X,\omega)$
such that $\|\Theta_n\|\leq c d_{k-p,n}\|\Theta\|$.
\end{lemma}
\begin{proof} We can assume that $\|\Theta\|=1$. The following constants $c$ and $c'$ 
depend only on $(X,\omega)$.
By \cite{DinhSibony3, DinhSibony4}, there exists
a sequence of smooth positive $(p,p)$-forms $\Psi_m$ which converges to
a current $\Psi\geq\Theta$ and such that $\|\Psi_m\|\leq c'$.
Hence the classes of $\Psi_m$
in $H^{p,p}(X,\C)$ are bounded. It follows that there exists $c>0$ such that 
$c\omega^p-\Psi_m$ is cohomologous to a positive closed form for every $m$.

Let $\Gamma_n$ be the graph of $f^n$ which
is an analytic subset of dimension $k$ of $X\times X$. 
If $\pi_1$, $\pi_2$ are the 
canonical projections of $X\times X$ onto its factors and $[\Gamma_n]$ is
the current of integration on $\Gamma_n$ then the positive closed 
currents 
$$(f^n)_*(\Psi_m):=(\pi_2)_*(\pi_1^*(\Psi_m)\wedge [\Gamma_n])$$
are well defined;
they have no mass on analytic sets because they have $\Lone$ coefficients. 
Since the mass of a positive closed current depends only on its cohomology class,
we have
\begin{eqnarray*}
\|(f^n)_*(\Psi_m)\| & \leq &  c\|(f^n)_*(\omega^p)\|
=c\int (f^n)_*(\omega^p)\wedge \omega^{k-p} \\
& = &  c \int \omega^p\wedge 
f^{n*}(\omega^{k-p})=cd_{k-p,n}.
\end{eqnarray*}
The previous integrals are computed on Zariski open sets where the forms are smooth.
Finally, we have 
$$\|\Theta_n\|\leq \limsup_{m\rightarrow\infty} 
\|(f^n)_*(\Psi_m)\|\leq cd_{k-p,n}.$$ 
\end{proof}

For every point $z$ out of an analytic subset of $X$, $f^{-1}(z)$ contains 
exactly $d_t$ points. We can define the Perron-Frobenius operator on the space 
of mesurable functions by
$$\Lambda\varphi(z):=d_t^{-1}f_*(\varphi)=
\frac{1}{d_t}\sum_{w\in f^{-1}(z)}\varphi(w).$$ 
Since $\nu_n$ gives no mass to analytic sets, the probability measure 
$d_t^{-n}f^{n*}(\nu_n)$ is well defined and we have 
$$\langle d_t^{-n}f^{n*}(\nu_n),\varphi\rangle = 
\langle\nu_n,\Lambda^n\varphi\rangle$$
for every $\varphi$ integrable with respect to $f^{n*}(\nu_n)$.

\begin{lemma} The operator $\Lambda: W_*^{1,2}\rightarrow W_*^{1,2}$ is well 
defined, bounded and continuous. 
When $\varphi$ is in $W^{1,2}_{**}$, so is $\Lambda\varphi$.
\end{lemma}
\begin{proof} Fix a Zariski open set $\Omega$ where $\Lambda\varphi$ is well defined.
The choice of $\Omega$ is not important. Fix also a compact set $K\subset\Omega$
big enough.
We have $\Lambda\varphi\in\Loneloc(\Omega)$ and $\|\Lambda\varphi\|_{\Lone(K)}\leq 
c\|\varphi\|_\Lone$, $c>0$. 

If $f_\star$ is the restriction of $f_*$ to a suitable Zariski open set, 
the Cauchy-Schwarz inequality implies that
$$i\partial(\Lambda\varphi)\wedge\overline\partial(\Lambda\varphi)=
d_t^{-2} f_\star(\partial\varphi)\wedge 
f_\star(\overline\partial\varphi)\leq
d_t^{-1}f_\star(i\partial\varphi\wedge \overline\partial\varphi) \leq d_t^{-1}
f_\star(\Theta)$$
where $\Theta$ is a positive closed current such that 
$i\partial\varphi\wedge\overline\partial\varphi\leq\Theta$
and $\|\Theta\|\leq\|\varphi\|^2_*$. 
Lemma 5.2 imply that $f_\star(\Theta)$ is positive closed and 
$\|f_\star(\Theta)\|\leq c\delta_1\|\varphi\|_*^2$. In particular,
$\d(\Lambda\varphi)$ has coefficients in $\Ltwo$. 
Proposition 3.1 implies that $\Lambda\varphi\in W^{1,2}$ and 
$$\|\Lambda\varphi\|_{\Lone}  \leq   
c'(\|\Lambda\varphi\|_{\Lone(K)}+ \|\d(\Lambda\varphi)\|_{\Ltwo})
 \leq  c''\|\varphi\|_*.$$
It follows that the operator $\Lambda: W_*^{1,2}\rightarrow W_*^{1,2}$
is bounded. The continuity follows.

When $\varphi\in W_{**}^{1,2}$, Lemma 4.7 implies that $\Lambda\varphi\in W_{**}^{1,2}$.
\end{proof}

Consider a function $\varphi\in W_*^{1,2}$. 
Define $c_0:=m(\varphi)$, $\varphi_0:=\varphi-c_0$, $c_{n+1}:=m(\Lambda\varphi_n)$
and $\varphi_{n+1}:=\Lambda\varphi_n -c_{n+1}$.
We have 
$$\Lambda^n\varphi=c_0+\cdots +c_n +\varphi_n.$$

\begin{lemma} There exists a constant $A>0$ independent of $\varphi$
such that $|c_n|\leq A\|\varphi\|_* \delta_{n-1}^{1/2} d_t^{-(n-1)/2}$ 
for $n\geq 1$ and 
$\|\varphi_n\|_*\leq A\|\varphi\|_*  \delta_n^{1/2} d_t^{-n/2}$
for $n\geq 0$.
\end{lemma}
\begin{proof} Let $\Theta$ be as above.
Define $\Theta_n:=d_t^{-n}(f^n)_\star(\Theta)$.
Lemma 5.2 implies that 
$\|\Theta_n\|\leq A^2\|\varphi\|_*^2\delta_n d_t^{-n}$ 
where $A>0$ is a constant. As in Lemma 5.3, we obtain 
$i\partial\varphi_n\wedge\overline\partial\varphi_n \leq \Theta_n$.
Since $m(\varphi_n)=0$, we have  $\|\varphi_n\|_*\leq \|\Theta_n\|^{1/2}$.
Hence $\|\varphi_n\|_*\leq A\|\varphi\|_*\delta_n^{1/2} d_t^{-n/2}$.
We then deduce from Lemma 5.3 that  
$|c_n|\leq\|\Lambda\varphi_{n-1}\|_*
\leq A\|\varphi\|_* \delta_{n-1}^{1/2} d_t^{-(n-1)/2}$
for some constant $A>0$.
\end{proof}

Define $c_\varphi:=\sum_{n\geq 0} c_n$. Lemmas 5.3 and 5.4 imply that 
$c_\varphi$ depends continuously on
$\varphi\in W_*^{1,2}$ and $|c_\varphi|\leq A\|\varphi\|_*$
for some constant $A>0$.
Now, assume that $\varphi$ is smooth. 
Then $\varphi_n\in W_{**}^{1,2}\cap \Linfty(X)$
with analytic singularities. The hypothesis on $\nu_n$ allows 
the following calculus 
$$\langle d_t^{-n}f^{n*}(\nu_n),\varphi \rangle = \langle \nu_n, 
\Lambda^n\varphi\rangle = c_0+\cdots +c_n +\langle \nu_n,\varphi_n \rangle$$
and Lemma 5.4 implies that $\langle \nu_n,\varphi_n\rangle$
tends to 0. Hence 
$\langle d_t^{-n}f^{n*}(\nu_n),\varphi \rangle$ converges 
to $c_\varphi$. It follows that $d_t^{-n}f^{n*}(\nu_n)$ converge weakly to a measure
$\mu$ independent of $(\nu_n)$. This measure $\mu$ is defined by
$$\langle \mu,\varphi \rangle := c_\varphi$$
for $\varphi$ smooth. The measure $\mu$ is WPC since 
we can extend it to a bounded continuous linear form on
$W_*^{1,2}$ by 
 $$\langle \mu,\varphi \rangle_* := c_\varphi.$$

\begin{lemma}
We have  $\langle \mu,\varphi \rangle_*=\langle \mu,\varphi \rangle$
for all $\varphi\in W_{**}^{1,2}$.
\end{lemma}
\begin{proof} 
We first prove that $\varphi$ is $\mu$-integrable. 
By Proposition 4.1, we have $\varphi=\varphi^+-\varphi^-$
with $\varphi^\pm\geq 0$, $\varphi^\pm\in W_{**}^{1,2}$ and
$\|\varphi^\pm\|_*\leq c\|\varphi\|_*$. Hence, we 
can assume that $\varphi$ is positive and continuous out of a compact set $Y$ of
W-capacity zero. If $K\subset X\setminus Y$ is a compact set, the integral of $\varphi$
on $K$ verifies
$$\langle \mu,\varphi \rangle_K \leq\limsup 
\langle d_t^{-n}f^{n*}(\omega^k),\varphi \rangle = \lim(c_0+\cdots+c_n) 
+   \limsup \langle \omega^k,\varphi_n \rangle
= c_\varphi.$$
Since $\mu(Y)=0$, $\varphi$ is $\mu$-integrable and $\langle\mu,\varphi\rangle\leq
c_\varphi\leq A\|\varphi\|_*$.

Define $b_n:=-\sum_{m>n}c_m$ and write 
$\varphi_n=\varphi_n^+-\varphi_n^-$ as in Proposition 4.1 with
$\varphi_n^\pm\geq 0$, $\varphi_n^\pm\in W_{**}^{1,2}$ 
and $\|\varphi_n^\pm\|\leq c\|\varphi_n\|_*$. 
We have $\langle\mu,\varphi_n^\pm\rangle \leq A\|\varphi_n^\pm\|_*
\leq Ac\|\varphi_n\|_*$.
Using the properties that $\mu$ is WPC and
$f^*(\mu)=d_t\mu$, we get
\begin{eqnarray*}
|\langle\mu,\varphi\rangle -c_\varphi| & = &  |\langle\mu,\Lambda^n\varphi
-c_\varphi\rangle| 
 \leq |b_n| +|\langle \mu, \varphi_n\rangle| \\
& \leq &  |b_n| + 
\langle \mu,\varphi_n^+\rangle + \langle \mu,\varphi_n^-\rangle\\
& \leq & |b_n| + 2Ac\|\varphi_n\|_*.
\end{eqnarray*}
Lemma 5.4 implies that the last expression tends to $0$. Hence
$\langle\mu,\varphi\rangle =c_\varphi$.
\end{proof}

\section{Mixing and central limit theorem}

In this section, we prove the geometric decay of correlations for the 
meromorphic maps that we are considering. We then deduce, via the 
Gordin-Liverani theorem, the CLT for Lipschitzian observables.

\begin{theorem} Let $f$ be as in Theorem 5.1. Then the equilibrium measure $\mu$
of $f$ is exponentially mixing in the following sense: there exists a
constant $A>0$ such that
$$|\langle\mu,(\psi\circ f^n) \varphi\rangle - \langle \mu,\psi\rangle
\langle \mu,\varphi\rangle |\leq A \|\psi\|_\infty\|\varphi\|_* \delta_n^{1/2} 
d_t^{-n/2}$$
for any $n\geq 0$ and 
every real-valued functions $\psi\in\Linfty(\mu)$ and 
$\varphi\in W_{**}^{1,2}$.
\end{theorem}
\begin{proof} Replacing $\varphi$ by $\varphi-c_\varphi$, we can assume that 
$c_\varphi=\langle \mu,\varphi\rangle =0$.
Assume also that $\|\varphi\|_*=1$ and $\|\psi\|_\infty=1$.  
The constants $A_i$ and $A$ below do not 
depend on $n$, $\psi$ and $\varphi$.

We keep the notations of the last section. 
Recall that $\|\varphi_n\|_*\leq A_1 \delta_n^{1/2} 
d_t^{-n/2}$ and  $\Lambda^n\varphi= b_n+\varphi_n$.
It is proved in \cite{DinhSibony3, DinhSibony4} 
that $((c\delta_n)^{1/n})$ is decreasing for
$c>0$ large enough. Hence, by Lemma 5.4, $|b_n|\leq A_2 \delta_n^{1/2} 
d_t^{-n/2}$.
As in Lemmas 5.4 and 5.5, we get 
$$\langle \mu, \varphi_n^\pm\rangle= \langle \mu, \varphi_n^\pm\rangle_*
\leq A_3\|\varphi_n^\pm\|_* \leq A_4 \delta_n^{1/2} d_t^{-n/2}.$$
We obtain from the relation $f^*(\mu)=d_t\mu$ that
\begin{eqnarray*}
|\langle\mu,(\psi\circ f^n) \varphi\rangle| & = & |\langle \mu, \psi \Lambda^n
\varphi\rangle| \\
&\leq & |b_n| \langle\mu,|\psi|\rangle + 
|\langle \mu, \psi \varphi_n\rangle|\\
& \leq & |b_n| + 
\langle \mu,|\psi|\varphi_n^+\rangle + \langle \mu,|\psi|\varphi_n^-\rangle\\
& \leq &  A_2 \delta_n^{1/2} 
d_t^{-n/2} + \langle \mu,\varphi_n^+\rangle + \langle \mu,\varphi_n^-\rangle\\
& \leq & A \delta_n^{1/2} d_t^{-n/2}.
\end{eqnarray*}
\end{proof}

\begin{corollary} Let $f$ and $\mu$ be as in Theorem 5.1. 
Let $\varphi$ be a bounded
function in $W_{**}^{1,2}$ such that $\langle\mu, \varphi\rangle=0$. 
Then if $\varphi$ is not a coboundary, it
satisfies the CLT. In particular, this holds for Lipschitzian functions $\varphi$
which satisfy $\langle \mu,\varphi\rangle=0$ and are not coboundaries. 
\end{corollary}
\begin{proof} If $\langle\mu,\varphi\rangle =0$, then 
$$\|\Lambda^n\varphi\|_{\Lone(\mu)}=
\sup_{\|\psi\|_\infty\leq 1}
|\langle \mu, \psi\Lambda^n\varphi\rangle|
=\sup_{\|\psi\|_\infty\leq 1}
|\langle \mu, (\psi\circ f^n)\varphi\rangle|.$$
Hence, the condition (1) is a direct consequence of Theorem 6.1. We then apply the 
Gordin-Liverani theorem.
\end{proof}

\begin{remark}\rm
{\bf (a)} All results still hold if we replace the closed current $\Theta$ in the
definition of $W^{1,2}_*$ by a $\partial\overline\partial$-closed current.
The regularization of $\partial\overline\partial$-closed currents
\cite{DinhSibony4} allows to prove an  analogue 
of Lemma 5.2.

{\bf (b)} We say that a positive measure $\nu$ is {\it PC} if its restriction 
to smooth real-valued functions can be extended to a linear continuous form on
$\DSH$. One can prove, 
using Demailly's regularization theorem as in Lemma 4.4, that d.s.h. functions
are $\nu$-integrable and the extension of $\nu$ 
is equal to the integral $\langle\nu, \varphi \rangle$ on $\varphi\in \DSH$. 
In \cite{DinhSibony2}, we proved that $\Lambda:\DSH\rightarrow\DSH$ 
is well defined, bounded and continuous. Using Hartogs lemma, one verifies that every 
$\varphi\in\DSH$ can be written as 
$\varphi=\varphi^+-\varphi^-$ with $\varphi^\pm\geq 0$
and $\|\varphi^\pm\|_\DSH\leq c\|\varphi\|_\DSH$, $c>0$. Then, we prove as above
that $\mu$ is PC and exponentially mixing in the sense that
$$|\langle\mu,(\psi\circ f^n) \varphi\rangle - \langle \mu,\psi\rangle
\langle \mu,\varphi\rangle |\leq A \|\psi\|_\infty\|\varphi\|_\DSH \delta_n 
d_t^{-n}$$
for any $n\geq 0$, $\psi\in\Linfty(\mu)$ and 
$\varphi\in \DSH$ (see \cite{DinhSibony2}). Hence, bounded d.s.h. functions 
satisfy the CLT.
One can construct bounded q.p.s.h. functions which are nowhere continuous. They
satisfy the CLT. 
As far as we know, this property is new even for rational fractions
in one variable.

{\bf (c)} The geometric decay of correlations was proved in \cite{DinhSibony1}
(see also   \cite{FornaessSibony})
for a large class of polynomial-like maps. This implies the CLT for 
bounded p.s.h. functions.
For rational fractions, the same problems were solved for H\"older
observables. See \cite{Denker, Haydn, Makarov, Ruelle} and the references therein.

{\bf (d)} In \cite{Leborgne}, Cantat and Leborgne have announced the CLT for holomorphic
endomorphisms of $\P^k$ and for H\"older continuous observables. 
Their idea is to use the
H\"older continuity of the Green function of $f$ and a
delicate estimate on modulus of annuli in order to apply the Gordin-Livenari theorem.
\end{remark}

\noindent
{\bf Acknowledgments.} 
The first named author thanks Guy David and 
Jean-Christophe L\'eger for their help during 
the preparation of this article.

\small

Tien-Cuong Dinh and Nessim Sibony, 
Math\'ematique - B\^at. 425, UMR 8628, \\
Universit\'e Paris-Sud, 91405 Orsay, France.\\ 
E-mails: (TienCuong.Dinh, Nessim.Sibony)@math.u-psud.fr
\end{document}